\documentclass[10pt]{article}
\usepackage{latexsym}
\usepackage{amsfonts}
\usepackage{enumerate}
\usepackage{multicol}
\usepackage{graphicx}
\usepackage{amssymb}
\usepackage{amsmath}
\usepackage{epic}

\title{On Asymptotic Approximate Groups of Integers  \\[.4in]}

\author{B\'{e}la Bajnok \\[.1in] Department of Mathematics, Gettysburg College \\
Gettysburg, PA 17325-1486 USA \\E-mail:  bbajnok@gettysburg.edu \\[.4in]}

\date{March 21, 2016}

\newtheorem{thm}{Theorem}

\newtheorem{lem}[thm]{Lemma}
\newtheorem{cor}[thm]{Corollary}
\newtheorem{prop}[thm]{Proposition}

\begin{document}

\maketitle

\begin{abstract}
 Let $r$ be a positive integer, and let $A$ be a nonempty finite set of at least two integers.  We let $\tilde{C}_r(A)$ denote the {\em asymptotic $r$-covering number} of $A$, that is, the smallest integer value of $l$ for which, for all sufficiently large positive integers $h$, the $rh$-fold sumset of $A$ is contained in at most $l$ translates of the $h$-fold sumset of $A$.  Nathanson proved that $\tilde{C}_r(A)$ is always at most $r+1$; here we extend this result to prove that $\tilde{C}_r(A)$ is always at least $r$, and determine all sets $A$ for which $\tilde{C}_r(A)=r$.

\end{abstract}


\thispagestyle{empty}

\section{Introduction}

We start by recalling some standard terms and notations.  Let $A$ be a nonempty finite set of integers.  If for some integers $a$ and $b$ with $a \leq b$, $A$ consists of all integers between $a$ and $b$, inclusive, we say that $A$ is an {\em interval} and write $A=[a,b]$; more generally, if for integers $a$, $d$, and $m$ with $d > 0$ and $m >0$, $A$ consists of all integers of the form $a+id$ with $i=0,1,\dots,m-1$, we say that $A$ is an {\em arithmetic progression} of common difference $d$, first term $a$, and size $m$.

For integers $x$ and $c$, we define the {\em translate of $A$ by $x$} as  $$A+x=\{a+x \mid a \in A\}$$ and the {\em $c$-fold   dilate of $A$} as $$c \cdot A=\{c \cdot a \mid a \in A\}.$$  For example, an arithmetic progression of common difference $d$, first term $a$, and size $m$ can be written as $d \cdot [0,m-1]+a$. 

We define the {\em sumset} of nonempty sets $A_1, A_2, \dots, A_h$ of integers by  $$A_1+A_2+ \cdots +A_h=\{a_1+a_2+ \cdots + a_h \mid a_i \in A_i \; \mbox{for} \; i=1,2,\dots,h\};$$  if $A_i=A$ for all $i=1,2,\dots,h$, we denote this sumset by $hA$ and call it the {\em $h$-fold sumset of $A$}.  (Note that $h \cdot A \subseteq hA$, but this is usually a proper containment.)

Let $r$ and $l$ be positive integers.  Nathanson \cite{Nat:2015a} defines $A$ to be an {\em $(r,l)$-approximate group}, if the $r$-fold sumset of $A$ is contained in the union of $l$ translates of $A$; that is, if there exists a set $X$ consisting of $l$ integers for which $rA \subseteq A+X$.  Note that any $(r,l)$-approximate group is an $(r,l')$-approximate group
for every integer $l' \geq l$.  It is, therefore, of interest to find the smallest integer $l$ for which $A$ is an $(r,l)$-approximate group; we here define this value as the {\em $r$-covering number of $A$} and denote it by $C_r(A)$.  Since $rA$ is a finite set, $C_r(A)$ is well-defined.

We make some trivial observations.  If $r=1$, then $rA=A$, so $A$ is a $(1,l)$-approximate group for any $l$; in particular, $C_1(A)=1$.  Also, when $|A|=1$, then $|rA|=1$, so $A$ is an $(r,l)$-approximate group for every $r$ and $l$, and $C_r(A)=1$.  For these reasons, below we assume that $r \geq 2$ and $|A| \geq 2$.

Let us consider two examples.

{\bf Example 1}  
Let $a_1$ and $a_2$ be any two distinct integers, and set $A=\{a_1,a_2\}$.  Then 
\begin{eqnarray*}
rA & = & \{(r-i)a_1+ia_2 \mid i=0,1,\dots,r\} \\ 
& \subseteq & A + \{(r-1-i)a_1+ia_2 \mid i=0,2,4,\dots, 2 \lceil (r-1)/2 \rceil \},
\end{eqnarray*} and so $C_r(A) \leq \lceil (r+1)/2 \rceil$.  (It is easy to see that equality holds here. In Section \ref{section Approximate groups of integers} we shall prove that $C_r(A) \geq \lceil (r+1)/2 \rceil$ for any set $A$.)

{\bf Example 2}    Let $a_1,a_2, \dots,a_m$ be positive integers so that $a_{i+1}>ra_i$ for $i=1,2,\dots,m-1$, and set
$A=\{a_1,a_2,\dots,a_m\}.$  Then $$|rA|={m+r-1 \choose r},$$ so we see that $A$ cannot be an $(r,l)$-approximate group for any $l < {m+r-1 \choose r} / m$.  In particular, we find that, for any pair of positive integers $r$ and $l$ (with $r \geq 2$), there are finite sets of integers $A$ with $C_r(A) >l$.

As these two examples demonstrate, the $r$-covering number of sets may vary widely.  We observe a much more subdued behavior if we consider sets `asymptotically'.  Namely, Nathanson  \cite{Nat:2015a} defines $A$ to be an {\em asymptotic $(r,l)$-approximate group}, if there is an integer $h_0$ so that $hA$ is an $(r,l)$-approximate group for every integer $h \geq h_0$.   The main result of \cite{Nat:2015a} is as follows.

\begin{thm} [Nathanson; cf.~\cite{Nat:2015a}] \label{Nathanson thm on asymp app group}
Every nonempty finite set of integers is an asymptotic $(r,r+1)$-approximate group.  
\end{thm}

We extend this result by establishing the following.

\begin{thm}  \label{my thm on asymp app group}
Let $A=\{a_1, a_2, \dots, a_m\}$ be a set of (at least two) integers so that $a_1 < a_2 \cdots < a_m$, and set $$d=\gcd(a_2-a_1, a_3-a_1, \dots, a_m-a_1).$$

\begin{enumerate}

\item If $A$ is an asymptotic $(r,l)$-approximate group, then $l \geq r$.

\item If $A$ is an asymptotic $(r,r)$-approximate group, then $a_2=a_1+d$ and $a_m=a_{m-1}+d$.

\end{enumerate}

\end{thm}

Analogously to the $r$-covering number $C_r(A)$ of $A$ introduced above, we define the {\em asymptotic $r$-covering number of $A$} as the smallest positive integer $l$ for which $A$ is an asymptotic $(r,l)$-approximate group; we denote this quantity by $\tilde{C}_r(A)$.  We may summarize Theorems \ref{Nathanson thm on asymp app group} and \ref{my thm on asymp app group} as follows.  

\begin{cor} \label{cor on asymp covering}
Let $A$ be any set of (at least two) integers.  Keeping the notations of Theorem \ref{my thm on asymp app group}, we have
$$\tilde{C}_r(A) = \left\{
\begin{tabular}{cl}
r & \mbox{if $a_2=a_1+d$ and $a_m=a_{m-1}+d$}; \\ \\
r+1 & \mbox{otherwise.}
\end{tabular}
\right.$$
\end{cor}

Contrasting our two examples above with Corollary \ref{cor on asymp covering} we can point out that, for any fixed $r$, the $r$-covering number $C_r(A)$ may take up infinitely many different values (depending on $A$), while its asymptotic version, $\tilde{C}_r(A)$, can only equal $r$ or $r+1$.

Approximate groups were introduced by Tao in \cite{Tao:2008a}, and have been investigated extensively since; see, for example, papers by Breuillard \cite{Bre:2014a};  Breuillard, Green, and Tao \cite{BreGreTao:2012a}; Green \cite{Gre:2012a}; Sanders \cite{San:2012a}; and their references.

\section{Approximate groups of integers} \label{section Approximate groups of integers}

We will use the following version of a well-known result.

\begin{lem} \label{lemma on min size of hA integers}
Let $h$ be a positive integer and $A$ be a set of integers.  If $A$ is an arithmetic progression, then $|hA|=hm-h+1$, otherwise $|hA|\geq hm$.

\end{lem}

{\em Proof:}  The first part of our claim is obvious.  The second part is obvious as well for $h=1$, and it is well-known for $h=2$: in fact, more generally, we know that, when $A_1$ and $A_2$ are nonempty finite sets of size at least two and at least one of them is not an arithmetic progression, then $|A_1+A_2| \geq |A_1|+|A_2|$ (see, for example, Lemma 1.3 and Theorem 1.3 in \cite{Nat:1996a}).  Suppose then that our claim holds for $A$ and some positive integer $k$, and consider $(k+1)A$.  Since both $A$ and $kA$ have size at least two, we have $$|(k+1)A|=|kA+A| \geq |kA|+|A| \geq km+m=(k+1)m,$$ thus our claim holds for all $h$ by induction.  $\Box$

\begin{prop} \label{app gr min l}
Let $r$ and $l$ be positive integers and $A$ be an $(r,l)$-approximate group of size $m$.  We then have $l \geq (rm-r+1)/m$; in fact, if $A$ is not an arithmetic progression, then $l \geq r$.

\end{prop}

{\em Proof:} Suppose that $X$ is a set of $l$ integers so that $rA \subseteq A+X$.  Then $$|rA| \leq |A+X| \leq |A| \cdot |X|=ml.$$ The result now follows from Lemma \ref{lemma on min size of hA integers}.  $\Box$

\begin{prop}  \label{app gr ap}
Let $r$ be a positive integer.  If $A$ is an arithmetic progression of size $m$, then $A$ is an $(r,l)$-approximate group for every $l \geq (rm-r+1)/m$; therefore, $$C_r(A)=\lceil (rm-r+1)/m \rceil.$$
In particular, we have $C_r(A) \leq r$ for any arithmetic progression $A$. 
\end{prop}

{\em Proof:}  Let $a$ and $d$ be integers such that $A=a+d \cdot [0,m-1]$.  Set $$X=(r-1)a+ md \cdot [0,l-1].$$  If $l \geq (rm-r+1)/m$, then 
$$(m-1)+(l-1)m \geq r(m-1),$$ and thus
\begin{eqnarray*}
rA & = & ra+d \cdot [0,r(m-1)] \\
& \subseteq & ra+d \cdot [0,(m-1)+(l-1)m] \\
& = & \left(a+d \cdot [0,m-1] \right) + \left((r-1)a+ md \cdot [0,l-1] \right) \\
& = & A+X,
\end{eqnarray*}
which implies the result. $\Box$

\begin{thm}
Let $r$ be an integer so that $r \geq 2$, and let $A$ be a finite set of integers of size $m \geq 2$.  Then $$C_r(A) \geq \lceil (r+1)/2 \rceil.$$  Furthermore, equality holds if, and only if
\begin{itemize}
  \item $m=2$ (and thus $A$ is an arithmetic progression), or
  \item $r=2$ and $A$ is an arithmetic progression (of any size), or 
  \item $r=4$ and $A$ is an arithmetic progression of size three.
\end{itemize} 

\end{thm}

{\em Proof:} If $A$ is an arithmetic progression, then $C_r(A)$ is given by Proposition \ref{app gr ap}; the verification of our claims is an easy exercise.

When $A$ is not an arithmetic progression, then by Proposition \ref{app gr min l} we have $C_r(A) \geq r \geq \lceil (r+1)/2 \rceil;$  it remains to be shown that at least one of the inequalities is strict.  If this were not the case, then we would have $r=2$, so we only need to verify that if $A$ is not an arithmetic progression, then it cannot be a $(2,2)$-approximate group.  This was done as part of Theorem 3 in \cite{Nat:1996a}; for the sake of completeness, we repeat the short proof here.

Suppose, indirectly, that $2A \subseteq A+X$ for some set $X$ of size two.  This implies that $|2A| \leq 2m$, so by Lemma \ref{lemma on min size of hA integers}, we must have $|2A| = 2m$ and thus $2A=A+X$.  But then $\min 2A=\min (A+X)$, where $\min 2A= 2 \cdot \min A$ and $\min (A+X)= \min A + \min X$.  Therefore, $\min A= \min X$; similarly, $\max A = \max X$.  This implies that $$\min A + \max X = \max A + \min X,$$ but then $|2A|=|A+X| \leq 2m-1$, a contradiction.  This completes our proof.  $\Box$

Regarding absolute upper bounds on $r$-covering numbers, we recall Example 2 to state that there are none.

\begin{prop}
Let $r$ be an integer so that $r \geq 2$.  Then for every positive integer $l$, there exists a nonempty finite set $A$ of integers so that  $C_r(A) \geq l$.  
\end{prop}

\section{Asymptotic arithmetic progressions of integers} \label{section Approximate arithmetic progressions of integers}

As we have seen, arithmetic progressions play a central role when studying approximate groups of integers.  In this section, as preparation for the study of asymptotic approximate groups, we investigate the `asymptotic' version of arithmetic sequences.  We call a subset $A$ of integers an {\em asymptotic arithmetic progression}, if $hA$ is an arithmetic progression for all sufficiently large $h$.

Our main result about asymptotic arithmetic progressions is as follows.

\begin{thm}  \label{thm on asymptotic arithmetic progressions}
Let $A=\{a_1, a_2, \dots, a_m\}$ be a set of (at least two) integers so that $a_1 < a_2 \cdots < a_m$, and set $$d=\gcd(a_2-a_1, a_3-a_1, \dots, a_m-a_1).$$ The following statements are equivalent.
\begin{enumerate}
  \item There is a positive integer $h_0$ for which $h_0A$ is an arithmetic progression.
  \item There is a positive integer $h_0$ for which $hA$ is an arithmetic progression for every integer $h \geq h_0$.
  \item $a_2=a_1+d$ and $a_m=a_{m-1}+d$.
\end{enumerate}

\end{thm}

Before our proof, we recall that each nonempty finite set of integers can be `normalized': the {\em normal form} of the set $A=\{a_1, a_2, \dots, a_m\}$ is the set $B=\{b_1=0,b_2,\dots,b_m\}$ defined by $b_i=(a_i-a_1)/d$ for $i=1,2,\dots,m$, where $$d=\gcd(a_2-a_1, a_3-a_1, \dots, a_m-a_1).$$  Note that $A=d \cdot B + a_1.$

\begin{lem} \label{lemma normalized for asymptotic arithmetic progressions}
Let $A$ be a nonempty finite set of integers, and let $B$ be its normal form.  Then for any positive integer $h$, $hA$ is an arithmetic progression if, and only if, $hB$ is an arithmetic progression.

\end{lem}

{\em Proof:}  Note that $hA$ is an arithmetic progression if, and only if, $hA-ha_1$ is an arithmetic progression, and $hB$ is an arithmetic progression if, and only if, $d \cdot hB$ is an arithmetic progression.  Since 
\begin{eqnarray*}
hA-ha_1 & = & \{\Sigma_{i=1}^m h_ia_i -ha_1 \mid \Sigma h_i=h\} \\
& = &   \{\Sigma_{i=1}^m (h_ia_i-h_ia_1) \mid \Sigma h_i=h\} \\
& = &  \{\Sigma_{i=1}^m d h_ib_i \mid \Sigma h_i=h\} \\
& = &  d \cdot \{\Sigma_{i=1}^m h_ib_i \mid \Sigma h_i=h\} \\
& = & d \cdot hB,
\end{eqnarray*} the result follows.  $\Box$

{\em Proof of Theorem \ref{thm on asymptotic arithmetic progressions}:}  We follow the scheme $1 \Rightarrow 3 \Rightarrow 2 \Rightarrow 1$; since $2 \Rightarrow 1$ is trivially true, we only need to establish the first two implications.

Suppose that $h_0A$ is an arithmetic progression for some positive integer $h_0$.  By Lemma \ref{lemma normalized for asymptotic arithmetic progressions}, we may assume that $A$ is in normal form, so we need to prove that $a_2=1$ and $a_m=a_{m-1}+1$.  

Note that $a_i \in h_0A$ for every $i=1,2,\dots,m$; in particular, the smallest element of $h_0A$ is $a_1=0$, and its smallest positive element is $a_2$.  Therefore, if $h_0A$ is an arithmetic progression, then each of its elements must be divisible by $a_2$.  In particular, $a_i$ is divisible by $a_2$ for every $i=1,2,\dots,m$.  Since $A$ is in normal form, this implies that $a_2=1$.  This means that $h_0A$ is an arithmetic progression with a common difference of 1.

Observe also that the largest element of $h_0A$ is $h_0a_m$, and the second-largest element is $(h_0-1)a_m+a_{m-1}$.  Since $h_0A$ is an arithmetic progression with a common difference of 1, we have $h_0a_m=(h_0-1)a_m+a_{m-1}+1$ and thus $a_m=a_{m-1}+1$.

Suppose now that we have a set $A$ with $a_2=a_1+d$ and $a_m=a_{m-1}+d$.  Then the normal form of $A$ is the set $B$ of the form $$B=\{0,1,b_3,\dots, b_{m-2},b-1,b\}$$ for some integers $1<b_3< \cdots < b_{m-2}<b-1.$ We shall prove that $$hB=\{0,1,2,\dots,hb\}$$  for every $h \geq b-2$; by Lemma \ref{lemma normalized for asymptotic arithmetic progressions}, this implies that $hA$ is an arithmetic progression for every integer $h \geq b-2$.  Since $\min hB=0$ and $\max hB=hb$, it suffices to prove that $[0,hb] \subseteq hB$.

Let $h \geq b-2$, let $0 \leq k \leq h$, and let $0 \leq j \leq k$. 
Note that the integer $$(k-j)\cdot (b-1)+j \cdot b$$ is an element of $kB$, so the interval $[k(b-1),kb]$ is contained in $kB$.  But $0 \in B$, so $kB \subseteq hB$,  and thus $[k(b-1),kb] \subseteq hB$.  

Now let $0 \leq t \leq h-k$.  Note that the integer $$(h-k-t) \cdot 0 + t \cdot 1 + k \cdot b$$ is an element of $hB$, and thus  the interval $[kb,kb+h-k]$ is contained in $hB$.  Combining this with the previous paragraph, we conclude that, for every integer $0 \leq k \leq h$, $hB$ contains the interval $[k(b-1),kb+h-k]$.  

But, since $h \geq b-2$, $$kb+h-k \geq (k+1)(b-1)-1,$$ we have $$[0,hb]=[0(b-1), hb+h-h]=\cup_{k=0}^h  [k(b-1),kb+h-k] \subseteq hB,$$ as claimed.  $\Box$

\section{Asymptotic approximate groups of integers} \label{section Asymptotic approximate groups of integers}

We start by proving that, in contrast to the $r$-covering number $C_r(A)$ of a set $A$, which may be less than $r$, the asymptotic $r$-covering number $\tilde{C}_r(A)$ of $A$ (of size at least two) is always at least $r$.

\begin{thm}
Let $A$ be any set of integers of finite size at least two, and let $r$ and $l$ be positive integers.  If $A$ is an asymptotic $(r,l)$-approximate group, then $l \geq r$. 

\end{thm}

{\em Proof:}  Suppose, indirectly, that $l<r$.  Let $m=|A| \geq 2$, $h$ be a positive integer, and suppose that $hA$ is an $(r,l)$-approximate group.  Let $X$ be a set of $l$ integers for which $rhA \subseteq hA+X$.  

We then have $$|rhA| \leq |hA| \cdot l;$$ since, by Lemma \ref{lemma on min size of hA integers}, $$|rhA|=|r (hA)| \geq r \cdot |hA| - r+1,$$
we get
$$|hA| \leq \frac{r-1}{r-l}.$$ 

But, again by Lemma \ref{lemma on min size of hA integers}, $|hA| \geq hm-h+1$, so $$hm-h+1 \leq \frac{r-1}{r-l},$$ or $$h \leq \frac{l-1}{(r-l)(m-1)}.$$  This implies that there can only be  finitely many values of $h$ for which $hA$ is an $(r,l)$-approximate group, which is a contradiction.  $\Box$

In the rest of this section, we characterize sets with asymptotic $r$-covering number equal to $r$.  We first reduce the problem to sets in normal form via the following lemma.

\begin{lem} \label{normal form approx group lemma}
Let $A$ be a nonempty finite set of integers, and let $B$ be its normal form.  Then, for any pair of positive integers $r$ and $l$, $A$ is an $(r,l)$-approximate group if, and only if, $B$ is an $(r,l)$-approximate group.

\end{lem}

{\em Proof:}  Suppose first that $B$ is an $(r,l)$-approximate group, and let $X$ be a set of $l$ integers for which $rB \subseteq B+X$.  Recall that, with $a_1=\min A$ and $$d=\gcd(a_2-a_1,\dots,a_m-a_1),$$ we have $A=d \cdot B + a_1$, and thus
$$rA=r(d \cdot B + a_1)=d \cdot rB + ra_1 \subseteq d \cdot (B+X) + ra_1= d \cdot B + d \cdot X + ra_1 = A + d \cdot X + (r-1)a_1.$$ 
Here $d \cdot X + (r-1)a_1$ is a set of $l$ integers, which implies that $A$ is an $(r,l)$-approximate group.

Conversely, assume that $A$ is an $(r,l)$-approximate group; let $A=\{a_1, a_2, \dots, a_m\}$ with $a_1$ and $d$ defined as above.  The covering number $C_r(A)$ of $A$ is at most $l$; let $X$ be a set of exactly $C_r(A)$ integers for which $rA \subseteq A+X$.

First, we prove that, for every $x \in X$, the integer $x-(r-1)a_1$ is divisible by $d$.  There must be an element $c \in rA$ for which $c \in A+x$, since otherwise we would have $$rA \subseteq A+(X \setminus \{x\}),$$  contradicting $|X|=C_r(A)$.  Suppose then that nonnegative integers $1 \leq j \leq m$ and $h_1, h_2, \dots, h_m$ are such that $\Sigma_{i=1}^m h_i =r$ and  
$$c=\Sigma_{i=1}^m h_ia_i =a_j+x.$$ 
Then
$$x-(r-1)a_1=\Sigma_{i=1}^m h_ia_i - ra_1 - (a_j - a_1) = \Sigma_{i=1}^m h_i(a_i - a_1) - (a_j - a_1),$$ which is then divisible by $d$ since $a_i-a_1$ is divisible by $d$ for every $i=1,2,\dots,m$.

Let $$Y=\{(x-(r-1)a_1)/d \mid x \in X\};$$ as we just proved, $Y$ is a set of $C_r(A)$ integers.  Furthermore,
$$d \cdot rB = r(d \cdot B+a_1) -ra_1 = rA -ra_1 \subseteq A+X -ra_1= d \cdot B +X -(r-1)a_1= d \cdot B + d \cdot Y,$$ so $rB \subseteq B+Y$.  Since $|Y|=C_r(A) \leq l$, we can conclude that $B$ is an $(r,l)$-approximate group, as claimed.  $\Box$

For the proof of Theorem \ref{characterization of asymptotic approx groups} below, we employ the following result, sometimes referred to as the Fundamental Theorem of Additive Number Theory.

\begin{thm} [Nathanson; cf.~Theorem 1.1 in \cite{Nat:1996a}] \label{Fundamental Theorem of Additive Number Theory}
Suppose that $A$ is a nonempty finite set of integers in normal form and $\max A=a_m$.  Then there are nonnegative integers $h_0$, $c$, and $d$, so that $hA$ contains the interval $[c,ha_m-d]$ for every integer $h \geq h_0$.

\end{thm}

For easier reference, we state our characterization of sets with $\tilde{C}_r(A)=r$ in the following form.

\begin{thm}  \label{characterization of asymptotic approx groups}
Let $A=\{a_1, a_2, \dots, a_m\}$ be a set of (at least two) integers so that $a_1 < a_2 \cdots < a_m$, and set $$d=\gcd(a_2-a_1, a_3-a_1, \dots, a_m-a_1).$$ Let $r$ be any positive integer.  The following statements are equivalent.
\begin{enumerate}
  \item $A$ is an asymptotic $(r,r)$-approximate group.
  \item $A$ is an asymptotic arithmetic progression.
  \item $a_2=a_1+d$ and $a_m=a_{m-1}+d$.
\end{enumerate}

\end{thm}
{\em Proof:}  We follow the scheme $1 \Rightarrow 3 \Rightarrow 2 \Rightarrow 1$.  The implication $3 \Rightarrow 2$ was already established in Theorem \ref{thm on asymptotic arithmetic progressions}.  Furthermore, if $A$ is an asymptotic arithmetic progression, that is, if $hA$ is an arithmetic progression for all sufficiently large integers $h$, then for these $h$  values $hA$ is an $(r,r)$-approximate group by Proposition \ref{app gr ap}, proving the implication $2 \Rightarrow 1$.  Therefore, we only need to prove that $1 \Rightarrow 3$.

Suppose then that $A$ is an asymptotic $(r,r)$-approximate group; by Lemma \ref{normal form approx group lemma} we may assume that $A$ is in normal form.  Our goal is to prove that $a_1=1$ and $a_m=a_{m-1}+1$.  We shall first prove that $a_2=1$.  

Suppose, indirectly, that $a_2 \geq 2$.  Then, by Theorem \ref{thm on asymptotic arithmetic progressions}, there is no positive integer $h$ for which $hA$ is an arithmetic progression, so, by Lemma \ref{lemma on min size of hA integers}, $|rhA| \geq r \cdot |hA|$ holds for every $h$.

Now let $h_0$, $c$, and $d$ be nonnegative integers for which, by Theorem \ref{Fundamental Theorem of Additive Number Theory}, $rhA$ contains the interval $[c,rha_m-d]$ for every integer $h \geq h_0/r$.  Then, with $$h_1=\max \{h_0/r, (c-1)/a_m, (d+1)/a_m\},$$ we have
$$[ha_m+1, (r-1)ha_m+1] \subseteq rhA$$ for every $h \geq h_1$.
  
Since $A$ is an asymptotic $(r,r)$-approximate group, we have an integer $h_2$ for which $hA$ is an $(r,r)$-approximate group for every integer $h \geq h_2$.  We set $$h_3=\max\{h_1,h_2,r/a_m\},$$ and choose $h$ to be any integer with $h \geq h_3$.

Let $X$ be a set of $r$ integers so that $rhA \subseteq hA+X$; we then have $|rhA| \leq |hA| \cdot r$.  Since above we have already established the reverse inequality, we must have $|rhA| = |hA| \cdot r$, and thus $rhA=hA+X$.  In particular, 
$$0=\min rhA = \min (hA+X) = \min hA + \min X = 0 + \min X,$$ so $\min X=0$; similarly, $\max X=(r-1)ha_m$.  Then we can write $X=\{x_1, x_2, x_3, \dots, x_{r-1}, x_r\}$, where
$$0=x_1<x_2<x_3 < \cdots < x_{r-1} < x_r=(r-1)ha_m.$$ 

Claim 1:  There is at least one value of $i \in \{1,2,\dots,r-1\}$ for which $x_{i+1}-x_i \geq ha_m$.

Proof of Claim 1:  This follows from the fact that $$\Sigma_{i=1}^{r-1} (x_{i+1}-x_i)= x_r-x_1=(r-1)ha_m.$$

Claim 2:  For all values of $i \in \{1,2,\dots,r-1\}$, we have $x_{i+1}-x_i \leq ha_m+1$.

Proof of Claim 2: Let us assume, to the contrary, that $x_{i+1}-x_i \geq ha_m+2$ for some $i$.  In this case the interval $$I=[x_i+ha_m+1,x_{i+1}-1]$$ is nonempty.  Note that $hA+X$ and $I$ are disjoint.  But each element of $I$ is at least $x_1+ha_m+1=ha_m+1$ and at most $x_r-1=(r-1)ha_m-1$, so 
$$I \subseteq [ha_m+1, (r-1)ha_m+1] \subseteq rhA,$$ which contradicts our assumption that $rhA = hA+X$.

Claim 3:  There is no value of $i \in \{1,2,\dots,r-1\}$ for which $x_{i+1}-x_i = ha_m$.

Proof of Claim 3: If, to the contrary, we were to have $x_{i+1}-x_i=ha_m$ for some $i$, then $x_{i+1}+1=x_i+ha_m+1$ would not be generated.  Since this integer is at least $ha_m+1$ and at most $(r-1)ha_m+1$, we get a contradiction as in the proof of Claim 2.

Claim 4: There is no value of $i \in \{1,2,\dots,r-1\}$ for which $x_{i+1}-x_i = 1$.

Proof of Claim 4:  If, to the contrary, we were to have $x_{i+1}-x_i=1$ for some value of $i$, then, by Claim 2, we would have
$$(r-1)ha_m=\Sigma_{i=1}^{r-1} (x_{i+1}-x_i) \leq 1+(r-2)(ha_m+1),$$ contradicting the fact that $h \geq h_3 \geq r/a_m$.  

Claim 5: $x_r-x_{r-1} \neq ha_m+1$.

Proof of Claim 5: Suppose, indirectly, that $x_r-x_{r-1} = ha_m+1$ and thus $x_{r-1}=(r-2)ha_m-1$.  Since $$(r-1)ha_m+1 \in rhA = hA+X,$$ we have $$a:=(r-1)ha_m+1-x_i \in hA$$ for some $i=1,2,\dots, r$.  We cannot have $i=r$, since then $a=1$, and $1 \not \in hA$ as the smallest positive element of $hA$ is $a_2 \geq 2$.  Also, we cannot have $i \leq r-1$, since then $$a \geq (r-1)ha_m+1-x_{r-1} = ha_m+2,$$ but the largest element of $hA$ is $ha_m$.  This proves our claim.

According to our five claims, we then must have a value of $i \in \{1,2,\dots,r-2\}$ for which $x_{i+1}-x_i = ha_m+1$ and $x_{i+2}-x_{i+1} \geq 2$.  Consider the integer $a=x_i+ha_m+2$.  We have $$ha_m+2 \leq a = x_{i+1}+1 \leq x_r = (r-1)ha_m,$$ so $a \in rhA =hA+X$.

However, we show that $a-X$ and $hA$ are disjoint.  Indeed, for all integers $0 \leq j \leq i$, we have $a-x_j \geq ha_m+2$; for all $i+2 \leq j \leq r$, we have $$a-x_j \leq x_i+ha_m+2-x_{i+2} \leq x_i+ha_m+2-(x_{i+1}+2)=-1;$$ and $a-x_{i+1}=1$.
But $$hA \subseteq \{0\} \cup [2,ha_m],$$ so  $a-X$ and $hA$ are disjoint, which is a contradiction.  This proves that $a_2=1$.

To prove that $a_m=a_{m-1}+1$, set $C=-A+a_m$.  Then $C$ is in normal form, and it is an easy exercise to verify that it is also an asymptotic $(r,r)$-approximate group; indeed, for any integer $h \geq h_3$, we have
\begin{eqnarray*}
rhC & = & rh(-A+a_m) \\
&=&-rhA +rha_m \\
&\subseteq &-(hA+X) +rha_m \\
&=& h(-A+a_m)+(-X+(r-1)ha_m) \\
&=& hC +(-X+(r-1)ha_m),
\end{eqnarray*} where $-X+(r-1)ha_m$ is a set of $r$ integers.  But then, by our argument above, the smallest positive element of $C$ equals 1; that is, $a_m-a_{m-1}=1$.  This completes our proof.
$\Box$

\end{document}